\documentclass[12pt, reqno]{amsart}
\usepackage{amsmath}
\usepackage{amssymb}
\usepackage{amsfonts}
\usepackage{amsthm}
\usepackage{calrsfs}
\usepackage{color}
\usepackage{hyperref}
\usepackage{graphicx}
\usepackage{amsaddr}
\usepackage{tensor}

\DeclareMathAlphabet{\pazocal}{OMS}{zplm}{m}{n}

\newtheorem{theorem}{Theorem}[section]
\newtheorem{proposition}[theorem]{Proposition}
\newtheorem{definition}[theorem]{Definition}

\newtheorem{example}[theorem]{Example}

\def\C{\pazocal{C}}
\def\N{\mathbb{N}}
\def\Z{\mathbb{Z}}
\def\P{\pazocal{P}}
\def\D{\pazocal{D}}

\def\Aut{\mathop{\mathrm{Aut}}}
\def\Atop{\mathop{\mathrm{Atop}}}
\def\Apar{\mathop{\mathrm{Apar}}}
\def\dev{\mathop{\mathrm{dev}}}
\def\Mult{\mathop{\mathrm{Mult}}}

\begin{document}

\title{Cubes of symmetric designs}

\author[V.~Kr\v{c}adinac, M.~O.~Pav\v{c}evi\'{c}, and K.~Tabak]{Vedran Kr\v{c}adinac$^1$, Mario Osvin Pav\v{c}evi\'{c}$^2$, and Kristijan Tabak$^3$}

\address{$^1$Faculty of Science, University of Zagreb, Bijeni\v{c}ka cesta~$30$, HR-$10000$ Zagreb, Croatia}

\address{$^2$Faculty of Electrical Engineering and Computing, University of Zagreb,
Unska~$3$, HR-$10000$ Zagreb, Croatia}

\address{$^3$Rochester Institute of Technology, Zagreb Campus,
D.~T.~Gavrana~$15$, HR-$10000$ Zagreb, Croatia}

\email{vedran.krcadinac@math.hr}
\email{mario.pavcevic@fer.hr}
\email{kristijan.tabak@croatia.rit.edu}

\thanks{This work has been supported by the Croatian Science Foundation
under the project $9752$.}



\date{September 10, 2023}

\maketitle

\begin{abstract}
We study $n$-dimensional matrices with $\{0,1\}$-entries ($n$-cubes) such
that all their $2$-dimensional slices are incidence matrices of symmetric designs.
A known construction of these objects obtained from difference sets is
generalized so that the resulting $n$-cubes may have inequivalent slices. For
suitable parameters, they can be transformed into $n$-dimensional
Hadamard matrices with this property. In contrast, previously known constructions
of $n$-dimensional designs all give examples with equivalent slices.
\end{abstract}

\section{Introduction}

In 1990, Warwick de Launey~\cite{dL90} introduced a framework for higher-dimen\-sional
combinatorial designs of various types. Among others, it encompasses symmetric
block designs, Hadamard matrices, and generalisations such as orthogonal designs
and weighing matrices. While $n$-dimensional Hadamard matrices have been
studied before \cite{PS71, PS79, JS80, HS81, HS80, SA81, SA84, YY86, YH87, YH88}
and after \cite{dL91, dLH93, KM02, dLS08, AAFR15} de~Launey's paper, and are even featured
in books \cite[\S 6]{SA85}, \cite[Chapter~5]{KH07}, \cite{YNX10}, \cite[Chapter~10]{ASEA11},
\cite[Chapter~11]{dLF11}, there seem to be no works dedicated to $n$-dimensional symmetric designs.

In this paper we study such objects under the name \emph{cubes of symmetric designs}.
These are $n$-dimensional $\{0,1\}$-matrices of order~$v$ such that all $2$-dimensional
slices are incidence matrices of $(v,k,\lambda)$ designs. Given an ordinary Hadamard
matrix of order~$v$, a proper $n$-dimensional Hadamard matrix of the same
order is obtained by the product contruction of Yang~\cite{YY86}. Thus, the spectrum of
orders~$v$ such that proper Hadamard matrices exist is the same for all
dimensions~$n\ge 2$. According to the famous Hadamard conjecture, it includes all orders
divisible by~$4$. No such construction is known for symmetric designs. There are parameters
such that symmetric designs exist, but $n$-dimensional cubes are not known
for $n\ge 3$, for example $(25,9,3)$. However, cubes of arbitrary dimension arise from
$(v,k,\lambda)$ difference sets, in analogy with the group development construction of
Hammer and Seberry~\cite{HS80}. We generalise this construction to commence from symmetric
designs that have difference sets as blocks, but are not necessarily developments. The
ensuing \emph{group cubes} have the interesting property that different slices may be
inequivalent $(v,k,\lambda)$ designs. In contrast, all known constructions of
proper $n$-dimensional Hadamard matrices and other types of $n$-dimensional designs
give examples with equivalent slices.

The outline of our paper is as follows. In Section~\ref{sec2} we recall the relevant definitions
and fix notation. Equivalence of cubes and their autotopy and autoparatopy groups are defined.
Invariants are introduced to be used  later on to distinguish inequivalent cubes.
Representations of cubes as orthogonal arrays and transversal designs are explained, enabling
the application of computational tools for incidence structures to cubes of symmetric designs.

Section~\ref{sec3} is devoted to the construction of $n$-cubes from difference sets. Properties
of such \emph{difference cubes} and their autotopy groups are proved. Equivalence of difference sets
and the corresponding difference cubes is compared. Distinctions are found between dimension
$n=2$ and dimensions $n\ge 3$.

The group cube construction is introduced in Section~\ref{sec4}. Examples not equivalent to
difference cubes are given for parameters $(21,5,1)$, $(16,6,2)$, and extended to
$(4^m,2^{m-1}(2^m-1),2^{m-1}(2^{m-1}-1))$ in Theorem~\ref{nondifcubes}. These are parameters
of Menon designs; thus, these cubes can be transformed into proper $n$-dimensional Hadamard
matrices with inequivalent slices.

Finally, Section~\ref{sec5} is devoted to computational results about cubes with small parameters.
A complete classification of group cubes is performed for $(16,6,2)$ and $(21,5,1)$.
Examples of non-difference group cubes are found for $(27,13,6)$, $(36,15,6)$, $(45,12,3)$,
$(63,31,15)$, $(64,28,$ $12)$, and $(96,20,4)$. More examples of $(16,6,2)$ cubes are constructed
by prescribing autotopy groups. Some of them have inequivalent parallel slices in all directions,
and are therefore not equivalent with any group cube.

Applications of $n$-dimensional Hadamard matrices include spectro\-scopy, error-correcting codes,
cryptography, and signal processing; see \cite[Section~5.3]{KH07} and \cite{YNX10}. Cubes of
symmetric designs could have similar applications, but this will be left for future work.

\section{Definitions}\label{sec2}

Let $V$ be a set of $v$ points. A \emph{$(v,k,\lambda)$ design}
over $V$ is a collection~$\D$ of $k$-subsets of $V$ called \emph{blocks}, such
that every pair of points is contained in exactly $\lambda$ blocks. The design
is \emph{symmetric} if the number of blocks is also~$v$. All designs
in this paper are symmetric. Given an ordering of the points $p_1,\ldots,p_v$ and
blocks $B_1,\ldots,B_v$, the \emph{incidence matrix} $A=(a_{ij})$ of the design is
defined by $a_{ij}=[p_i\in B_j]$. Here and in the sequel, we use the Iverson bracket
$[P]$ which takes the value~$1$ if~$P$ is true and~$0$ otherwise; see~\cite{DK92}.
It is known that a $v\times v$ matrix~$A$ with $\{0,1\}$-entries is the incidence
matrix of a $(v,k,\lambda)$ design if and only if $A\cdot A^t=(k-\lambda)I+\lambda J$
holds, where $I$ is the identity matrix, and $J$ is the all-ones matrix. Furthermore,
the transposed matrix $A^t$ satisfies the same equation and is an incidence matrix of
the dual design $\D^t$. We refer to~\cite{EL83} for these and other
results about symmetric designs.

Let $\{1,\ldots,v\}^n$ be the Cartesian product of $n$ copies of $\{1,\ldots,v\}$.
An $n$-dimensional \emph{incidence cube} of order $v$ is a function $C:\{1,\ldots,v\}^n\to
\{0,1\}$. We will talk about $n$-cubes of order $v$ for short. Thus, a $2$-cube
is a $v\times v$ matrix with $\{0,1\}$-entries, a $3$-cube is a $3$-dimensional
array of zeros and ones, etc. Given a pair of distinct integers $(x,y)\in\{1,\ldots,n\}^2$,
a \emph{slice} of the $n$-cube $C$ is the matrix obtained by varying the coordinates in
positions $x$ and $y$, and taking some fixed values
$i_1,\ldots,i_{x-1},i_{x+1},\ldots,i_{y-1},i_{y+1},\ldots,i_n\in \{1,\ldots,v\}$
for the remaining coordinates. In other words, it is the restriction of
$C$ to the set $$\{i_1\}\times\cdots\times \{i_{x-1}\}\times V \times \{i_{x+1}\}\times
\cdots \times \{i_{y-1}\}\times V\times \{i_{y+1}\}\times \cdots \times \{i_n\}.$$
We allow $x>y$, in which case the order of the factors indexed by $x$ and $y$ should be
reversed. The slices corresponding to $(x,y)$ and $(y,x)$ are transposed matrices.

\begin{definition}\label{maindef}
An \emph{$n$-dimensional cube of symmetric $(v,k,\lambda)$ designs}
is an $n$-cube of order $v$ such that all of its slices are incidence
matrices of $(v,k,\lambda)$ designs. The set of all such $n$-cubes
will be denoted by $\C^n(v,k,\lambda)$.
\end{definition}

This is a special case of de Launey's proper $n$-dimensional transposable designs
$(v,\Pi_R,\Pi_C,\beta,S)^n$, see \cite[Definitions 2.1, 2.6, 2.7, and Example 2.2]{dL90}.
Here is a concrete example of such an object.

\begin{example}\label{motivation}
Consider the following incidence matrix of the Fano plane, i.e.\ the  $(7,3,1)$ design:
$$A_1 = \left(\begin{array}{ccccccc}
   1 & 1 & 0 & 1 & 0 & 0 & 0 \\
   1 & 0 & 1 & 0 & 0 & 0 & 1 \\
   0 & 1 & 0 & 0 & 0 & 1 & 1 \\
   1 & 0 & 0 & 0 & 1 & 1 & 0 \\
   0 & 0 & 0 & 1 & 1 & 0 & 1 \\
   0 & 0 & 1 & 1 & 0 & 1 & 0 \\
   0 & 1 & 1 & 0 & 1 & 0 & 0  \\
\end{array}\right).$$
Let $A_2$ be the matrix obtained by a cyclic shift upwards of rows of $A_1$.
This is another incidence matrix of the Fano plane. Continue in the same
way to get incidence matrices $A_3$ to $A_7$:
$$A_2=\left(\begin{array}{ccccccc}
    1 & 0 & 1 & 0 & 0 & 0 & 1 \\
    0 & 1 & 0 & 0 & 0 & 1 & 1 \\
    1 & 0 & 0 & 0 & 1 & 1 & 0 \\
    0 & 0 & 0 & 1 & 1 & 0 & 1 \\
    0 & 0 & 1 & 1 & 0 & 1 & 0 \\
    0 & 1 & 1 & 0 & 1 & 0 & 0 \\
    1 & 1 & 0 & 1 & 0 & 0 & 0  \\
\end{array}\right), \quad
A_3=\left(\begin{array}{ccccccc}
    0 & 1 & 0 & 0 & 0 & 1 & 1 \\
    1 & 0 & 0 & 0 & 1 & 1 & 0 \\
    0 & 0 & 0 & 1 & 1 & 0 & 1 \\
    0 & 0 & 1 & 1 & 0 & 1 & 0 \\
    0 & 1 & 1 & 0 & 1 & 0 & 0 \\
    1 & 1 & 0 & 1 & 0 & 0 & 0 \\
    1 & 0 & 1 & 0 & 0 & 0 & 1  \\
\end{array}\right),$$
$$A_4=\left(\begin{array}{ccccccc}
    1 & 0 & 0 & 0 & 1 & 1 & 0 \\
    0 & 0 & 0 & 1 & 1 & 0 & 1 \\
    0 & 0 & 1 & 1 & 0 & 1 & 0 \\
    0 & 1 & 1 & 0 & 1 & 0 & 0 \\
    1 & 1 & 0 & 1 & 0 & 0 & 0 \\
    1 & 0 & 1 & 0 & 0 & 0 & 1 \\
    0 & 1 & 0 & 0 & 0 & 1 & 1  \\
\end{array}\right), \quad
A_5=\left(\begin{array}{ccccccc}
    0 & 0 & 0 & 1 & 1 & 0 & 1 \\
    0 & 0 & 1 & 1 & 0 & 1 & 0 \\
    0 & 1 & 1 & 0 & 1 & 0 & 0 \\
    1 & 1 & 0 & 1 & 0 & 0 & 0 \\
    1 & 0 & 1 & 0 & 0 & 0 & 1 \\
    0 & 1 & 0 & 0 & 0 & 1 & 1 \\
    1 & 0 & 0 & 0 & 1 & 1 & 0  \\
\end{array}\right),$$
$$A_6=\left(\begin{array}{ccccccc}
    0 & 0 & 1 & 1 & 0 & 1 & 0 \\
    0 & 1 & 1 & 0 & 1 & 0 & 0 \\
    1 & 1 & 0 & 1 & 0 & 0 & 0 \\
    1 & 0 & 1 & 0 & 0 & 0 & 1 \\
    0 & 1 & 0 & 0 & 0 & 1 & 1 \\
    1 & 0 & 0 & 0 & 1 & 1 & 0 \\
    0 & 0 & 0 & 1 & 1 & 0 & 1  \\
\end{array}\right), \quad
A_7=\left(\begin{array}{ccccccc}
    0 & 1 & 1 & 0 & 1 & 0 & 0 \\
    1 & 1 & 0 & 1 & 0 & 0 & 0 \\
    1 & 0 & 1 & 0 & 0 & 0 & 1 \\
    0 & 1 & 0 & 0 & 0 & 1 & 1 \\
    1 & 0 & 0 & 0 & 1 & 1 & 0 \\
    0 & 0 & 0 & 1 & 1 & 0 & 1 \\
    0 & 0 & 1 & 1 & 0 & 1 & 0 \\
\end{array}\right).$$
Now put all of these matrices on top of each other: first $A_1$, then $A_2$ and end up with $A_7$ as the
top layer. We get a $3$-cube of order $7$ such that the horizontal slices are incidence matrices of the $(7,3,1)$
design. Because of the cyclic nature of the construction, this also holds for slices in the other two directions.
If we look at a particular row of $A_1,\ldots,A_7$, say row $j$, we get the matrix~$A_j$. The matrices are symmetric
and the conclusion also holds if we look at a particular column.
\end{example}

We now define equivalence of cubes of symmetric designs; cf.~\cite{KM02, dLS08} for
the corresponding concepts for $n$-dimensional Hadamard matrices. The direct product
of symmetric groups $(S_v)^n=S_v\times\ldots\times S_v$ acts
on the set $\C^n(v,k,\lambda)$ by permuting indices: for $\alpha=(\alpha_1,\ldots,\alpha_n)\in (S_v)^n$,
$$C^\alpha (i_1,\ldots,i_n) = C(\alpha_1^{-1}(i_1),\ldots,\alpha_n^{-1}(i_n)).$$
The orbits of this action
are the \emph{isotopy classes} of cubes. The stabiliser of a cube
$C\in \C^n(v,k,\lambda)$ is its \emph{autotopy group} $\Atop(C)$. Furthermore,
permutations $\gamma\in S_n$ act by \emph{conjugation}, i.e.\ changing the order of the indices:
$$C^\gamma (i_1,\ldots,i_n) = C(i_{\gamma^{-1}(1)},\ldots,i_{\gamma^{-1}(n)}).$$
A cube such that $C^\gamma=C$ for all $\gamma\in S_n$
is \emph{totally symmetric}. The cube from Example~\ref{motivation} has this property.
The combination of isotopy and conjugation is called \emph{paratopy} and is the natural
action of the wreath product $S_v\wr S_n$ on $\C^n(v,k,\lambda)$. Two cubes are considered
\emph{equivalent} if they can be mapped onto each other by paratopy. The corresponding
classes are the \emph{main classes} and the stabiliser of a cube $C$ in $S_v\wr S_n$ is
its \emph{autoparatopy group} $\Apar(C)$.

The terminology is borrowed from Latin squares~\cite{KD15}. Indeed, a Latin
square $L=(\ell_{i_1 i_2})$ of order $v$ is equivalent to a $3$-cube
$C\in \C^3(v,1,0)$ by $C(i_1,i_2,i_3)=[\ell_{i_1 i_2}=i_3]$,
since incidence matrices of symmetric designs with $k=1$ are just permutation
matrices. The equivalence is easily generalised to higher dimensions, but notice
that different authors have considered different concepts under the name
``Latin hypercubes''. The appropriate definition for our purposes is the one
used in~\cite{MW08}; see the discussion and references therein for other
related concepts. In the case $k=2$, symmetric designs are possible only for
parameters $(3,2,1)$. There is a unique $(3,2,1)$ design up to isomorphism and
it is cyclic. By Theorem~\ref{tmdifcube}, cubes of $(3,2,1)$ designs exist for
all dimensions $n\ge 2$. In the sequel we consider only designs with $k\ge 3$.

To distinguish inequivalent cubes, we use invariants based on parallel
slices. For dimension $n=2$, cubes are simply incidence matrices of symmetric
designs and equivalence means that the designs are isomorphic or dually
isomorphic. For $n>2$, we can choose the varying coordinates in ${n\choose 2}$
ways and substitute the fixed coordinates in $v^{n-2}$ ways. Slices are \emph{parallel}
if all but one of the fixed coordinates agree. We consider $v$ mutually
parallel slices as a multiset of (dual) isomorphism types of $(v,k,\lambda)$ designs.
The multiset of ${n\choose 2}(n-2)v^{n-3}$ such multisets must agree for equivalent
cubes, i.e.\ this is a paratopy invariant. We may also look at weaker invariants,
for example orders of the full automorphism groups of parallel slices. In either
case, the multiset of such multisets agreeing is only a necessary, but not sufficient
condition for equivalence of cubes; see the comment after Example~\ref{pg24}.

An incidence cube $C:\{1,\ldots,v\}^n\to\{0,1\}$ is the characteristic function
of a set of $n$-tuples
$$\overline{C}=\{ (i_1,\ldots,i_n) \in \{1,\ldots,v\}^n \mid C(i_1,\ldots,i_n)=1\}.$$
This is a one-to-one correspondence. If $C\in \C^n(v,k,\lambda)$, then $\overline{C}$
is an orthogonal array with parameters $OA(kv^{n-1},n,v,n-1)$ of index~$k$; see~\cite{HSS99}
or~\cite[Section III.6]{CD07} for the definition. It has the additional property that for any
choice of coordinate positions $x,y\in\{1,\ldots,n\}$, $x\neq y$, and values
$i_1,\ldots,i_{x-1},i_{x+1},\ldots,i_{y-1},i'_y,i''_y,i_{y+1},\ldots,i_n \in \{1,\ldots,v\}$,
$i'_y\neq i''_y$, there are exactly~$\lambda$ values $i\in \{1,\ldots,v\}$ such that
$(i_1,\ldots,i_{x-1},$ $i,i_{x+1},\ldots,i'_{y},\ldots,i_n)\in \overline{C}$ and
$(i_1,\ldots,i_{x-1},i,i_{x+1},\ldots,i''_{y},\ldots,i_n)\in \overline{C}$. Every orthogonal
array with this property and the parameters above corresponds to an $n$-dimensional
cube of $(v,k,\lambda)$ designs. Thus, we can represent cubes of symmetric designs by
orthogonal arrays. In practice we use an equivalent representation by transversal
designs: $\overline{C}$ is represented as a collection of $n$-subsets
$\{i_1,v+i_2,2v+i_3,\ldots,(n-1)v+i_n\}$ instead of $n$-tuples $(i_1,\ldots,i_n)$.
Now we have an incidence structure of $nv$ points and $kv^{n-1}$
blocks such that the usual notion of isomorphism agrees with paratopy of cubes.
This has the advantage that we can check equivalence of cubes and compute
autoparatopy groups by tools available in the computer algebra system GAP~\cite{GAP4},
notably the graph isomorphism programs \emph{nauty} and \emph{Traces}~\cite{MP14}.
By coloring the $n$ classes of points $\{1,\ldots,v\},\ldots,\{(n-1)v+1,\ldots,nv\}$
we can also check isotopy of cubes, compute autotopy groups, and use the Kramer-Mesner
approach~\cite{KM76} to construct cubes with prescribed autotopy groups. These and
other tools for handling cubes of symmetric designs are available in our GAP package
\emph{Prescribed Automorphism Groups}~\cite{PAG}.

\section{Difference cubes}\label{sec3}

Let $G$ be a multiplicatively written group of order~$v$. A $(v,k,\lambda)$
\emph{difference set} in $G$ is a $k$-subset $D\subseteq G$ such that every
$g\in G$, $g\neq 1$ can be written as a left difference $g=d_1^{-1}d_2$ for
exactly $\lambda$ pairs $(d_1,d_2)\in D\times D$. This is equivalent with the
corresponding property for right differences $g=d_1 d_2^{-1}$~\cite{RHB55}.
The \emph{development} of~$D$, i.e.\ the set of left translates
$\dev D=\{gD \mid g\in G\}$ is a $(v,k,\lambda)$ design over~$G$.
The right translates $\{Dg\mid g\in G\}$ also form a design
isomorphic to the dual $(\dev D)^t$. Difference sets give rise to cubes
of symmetric designs of arbitrary dimension. The following theorem is a
special case of~\cite[Theorem~2.9]{dL90}. An analogous theorem was originally
proved for group developed $n$-dimensional proper Hadamard
matrices~\cite[Theorem~4]{HS80}.

\begin{theorem}\label{tmdifcube}
Let $D$ be a $(v,k,\lambda)$ difference set in the group~$G$.
Order the group elements as $g_1,\ldots,g_v$. Then the function
\begin{equation}\label{dcube}
C(i_1,\ldots,i_n)=[g_{i_1}\cdots g_{i_{n}} \in D]
\end{equation}
is an $n$-dimensional cube of $(v,k,\lambda)$ designs.
\end{theorem}

\begin{proof}
First consider the case $n=2$. Notice that the design $\dev D$ with points
ordered as $g_1,\ldots,g_v$ and blocks ordered as $g_1D,\ldots,g_vD$ has
incidence matrix $A=(a_{ij})$ given by $a_{ij}=[g_i \in g_j D]$.
The matrix $C(i,j)=[g_i g_j \in D]$ is an incidence matrix of the dual
design $(\dev D)^t$ with a different ordering of columns, since
$g_ig_j\in D \iff g_j\in g_i^{-1}D$ holds. In the general case,
a slice of the $n$-cube~\eqref{dcube} is a matrix of the form
$a_{ij}=[h_1 g_i h_2 g_j h_3 \in D]$ or $a_{ij}=[h_1 g_j h_2 g_i h_3 \in D]$
for some fixed $h_1,h_2,h_3\in G$. This is just a permutation of
rows and columns of the incidence matrix of $(\dev D)^t$ or $\dev D$,
and therefore~\eqref{dcube} is a $(v,k,\lambda)$ cube of dimension~$n$.
\end{proof}

The general construction of Theorem~2.9 in~\cite{dL90} relies on
\emph{collapsable functions}, which were later specialised to
\emph{abelian extension functions} and \emph{$2$-cocycles}~\cite{dLH93}.
Cocyclic Hadamard matrices are an important and widely studied class~\cite{KH07},
more general than group developed Hadamard matrices. The latter are
necessarily regular, and therefore of order $v=4u^2$ for some $u\in \N$.
There is no such restriction for cocyclic Hadamard matrices, which are conjectured to
exist for all orders divisible by~$4$ \cite[Conjecture~3.6]{dLH93}. Regarding
symmetric designs, the symbols $0$ and $1$ in their incidence matrices cannot
be exchanged and there are no nontrivial cocycles. In this case
Theorem~2.9 of~\cite{dL90} reduces to the difference set construction of
Theorem~\ref{tmdifcube}.

We shall call cubes arising from Theorem~\ref{tmdifcube} \emph{difference cubes}.
A \emph{non-difference cube} is a cube of symmetric designs not equivalent to any
difference cube. The following properties of difference cubes are easy consequences
of the definition.

\begin{proposition}\label{cordc}
All slices of an $n$-dimensional difference cube are incidence
matrices of the same $(v,k,\lambda)$ design up to
isomorphism and duality.
\end{proposition}

\begin{proposition}
A difference cube constructed from an abelian group is totally symmetric,
i.e.\ invariant under any conjugation.
\end{proposition}

A design $\dev D$ constructed from a difference set $D\subseteq G$ has $G$ as an automorphism
group acting regularly on the points and blocks. An $n$-dimensional difference cube has
the direct product $G^{n-1}=G\times\ldots\times G$ as autotopy group acting as follows.
For $a,b\in G$, let $\tensor*[_a]{\alpha}{_b}\in S_v$ be the permutation defined
by \hbox{$\tensor*[_a]{\alpha}{_b}(i)=j \iff ag_i b^{-1}=g_j$}. Given $(a_1,\ldots,a_{n-1})\in G^{n-1}$,
let $a_0=a_n=1$ and $\alpha_i=\tensor*[_{a_{i-1}}]{\alpha}{_{a_i}}$, for $i=1,\ldots,n$.
Then $(a_1,\ldots,a_{n-1})\mapsto (\alpha_1,\ldots,\alpha_n)=:\alpha$ is a group embedding
$G^{n-1}\hookrightarrow (S_v)^n$ and its image is an autotopy group
of the difference cube~\eqref{dcube}:
\begin{align*}
C(i_1,\ldots,i_n) &= [g_{i_1}\cdots g_{i_{n}} \in D] \\[2mm]
 &= [g_{i_1}a_1\cdot a_1^{-1}g_{i_2}a_2\cdot\ldots\cdot a_{n-2}^{-1}g_{i_{n-1}}a_{n-1}\cdot a_{n-1}^{-1}g_{i_{n}} \in D]\\[2mm]
 &= [g_{\alpha_1^{-1}(i_1)}\cdots g_{\alpha_n^{-1}(i_{n})} \in D] =C^\alpha (i_1,\ldots,i_n).
\end{align*}

A \emph{multiplier} of the difference set $D\subseteq G$ is a group automorphism $\varphi:G\to G$
such that $\varphi(D)=aD$ for some $a\in G$. The set of all multipliers of $D$ is a subgroup
$\Mult(D)\le \Aut(G)$ acting on the design $\dev D$ as automorphism group. This
group also acts on the difference cube~\eqref{dcube} as autotopy group:
\begin{align*}
g_{i_1}\cdots g_{i_{n}} \in D &\iff \varphi(g_{i_1}\cdots g_{i_{n}})\in \varphi(D) \iff
\varphi(g_{i_1})\cdots \varphi(g_{i_{n}})\in aD\\[2mm]
&\iff a^{-1}\varphi(g_{i_1})\cdots \varphi(g_{i_{n}})\in D.
\end{align*}
In the last line we see how to permute the indices to leave the difference cube~\eqref{dcube}
invariant. The combined actions of $G^{n-1}$ and $\Mult(D)$ give an autotopy group of~\eqref{dcube}
isomorphic to their semidirect product. This proves the following theorem.

\begin{theorem}
The difference cube~\eqref{dcube} has an autotopy group isomorphic to a semidirect
product $(G^{n-1})\rtimes \Mult(D)$.
\end{theorem}

Two difference sets $D_1\subseteq G_1$ and $D_2\subseteq G_2$ are \emph{equivalent}
if there is a group isomorphism $\varphi:G_1\to G_2$ such that $\varphi(D_1)=aD_2$
for some $a\in G_2$. Equivalent difference sets clearly give rise to isomorphic
designs and isotopic difference cubes. For designs, the converse does not hold:
$\dev D_1$ and $\dev D_2$ can be isomorphic even if $D_1$ and $D_2$ are not equivalent.
For example, according to the GAP package \emph{DifSets}~\cite{DP19} there are $27$
inequivalent $(16,6,2)$ difference sets in $12$ of the $14$ groups of order~$16$
(see Table~\ref{gc16}). On the other hand, there are only three $(16,6,2)$ designs
up to isomorphism and duality~\cite{CD07}. By a computation in PAG~\cite{PAG},
the $27$ difference $3$-cubes $(16,6,2)$ are not equivalent (paratopic).
The conclusion holds for all dimensions $n\ge 3$ because any sub-$3$-cube of a
difference $n$-cube is equivalent with the corresponding difference $3$-cube.

The smallest examples of equivalent difference $n$-cubes ($n\ge 3$) obtained from inequivalent
difference sets have parameters $(27,13,6)$. Difference sets exist in two of the five groups of
order~$27$. There is a unique $(27,13,6)$ difference set in $\Z_3\times\Z_3\times \Z_3$
and two inequivalent difference sets in $\Z_9\rtimes \Z_3$~\cite{DP19}. Another computation in
PAG~\cite{PAG} shows that the difference $3$- and $4$-cubes constructed from the two difference
sets in $\Z_9\rtimes \Z_3$ are equivalent, but not isotopic. See Table~\ref{tabndgc} for other
examples of equivalent difference cubes coming from inequivalent difference sets. For dimensions
$n\ge 3$, we did not find a single example of isotopic difference cubes coming from inequivalent
difference sets, nor of equivalent difference cubes coming from difference sets in nonisomorphic
groups. Both situations are possible for $n=2$.

\section{Group cubes}\label{sec4}

In this section we generalise the construction of difference cubes of
Theorem~\ref{tmdifcube}. Let $G=\{g_1,\ldots,g_v\}$ be a group of order~$v$
and $\D=\{B_1,\ldots,B_v\}$ a $(v,k,\lambda)$ design over~$G$ such that all
blocks $B_i$ are $(v,k,\lambda)$ difference sets. The design $\D$ could be the
development of a difference set, say $B_i=g_i^{-1} B_1$, $i=1,\ldots,v$. In this
case the corresponding difference cube can be written as $C(i_1,\ldots,i_n)=[g_{i_1}\cdots g_{i_{n}} \in B_1]=$
$[g_{i_2}\cdots g_{i_{n}} \in g_{i_1}^{-1}B_1] = [g_{i_2}\cdots g_{i_{n}} \in B_{i_1}]$.
The last formula gives a $(v,k,\lambda)$ cube even if $\D$ is not a development.

\begin{theorem}\label{tmgrcube}
Let $G=\{g_1,\ldots,g_v\}$ be a group and $\D=\{B_1,\ldots,B_v\}$ a $(v,k,\lambda)$
design with all of its blocks being $(v,k,\lambda)$ difference sets in~$G$. Then
\begin{equation}\label{gcube}
C(i_1,\ldots,i_n)=[g_{i_2}\cdots g_{i_{n}} \in B_{i_1}]
\end{equation}
is an $n$-dimensional cube of $(v,k,\lambda)$ designs.
\end{theorem}

\begin{proof}
If the index $i_1$ is fixed, then~\eqref{gcube} is just the $(n-1)$-dimensional
difference cube constructed from the difference set $B_{i_1}$, so any such slice
is an incidence matrix of its development. On the other hand, if $i_1$ and one
of the remaining indices are varied, then the corresponding slice is an incidence
matrix of the design $\D^t$ with reordered columns.
\end{proof}

We shall call cubes constructed from Theorem~\ref{tmgrcube} \emph{group cubes}.
The main question is whether there exist designs that are not developments,
but all of their blocks are difference sets. In this case the construction may give
non-difference cubes.
Our first example is for parameters $(21,5,1)$ of the
projective plane of order~$4$.

\begin{example}\label{pg24}
There are two groups of order $21$: the Frobenius group
$F_{21}=\langle a,b \mid a^3=b^7=1,\, ba=ab^2\rangle$ and the cyclic group~$\Z_{21}$.
Both groups allow one $(21,5,1)$ difference set up to equivalence~\cite{DP19}.
We denote the corresponding $3$-cubes $C_1$ and $C_2$, respectively. Using PAG~\cite{PAG},
we computed $|\Atop(C_1)|=1323$ and $|\Atop(C_2)|=2646$. Thus, $C_1$ and $C_2$ are not
equivalent, and neither are the corresponding $n$-cubes for $n>3$ because all of their
sub-$3$-cubes are equivalent either with $C_1$ or with $C_2$. This shows that there are
exactly two inequivalent difference cubes in $\C^n(21,5,1)$ for every $n\ge 3$.

Here is a $(21,5,1)$ design over $F_{21}$ with all blocks being difference sets,
which is not the development of any of its blocks:
$$\begin{array}{l}
\{\,\,  \{ 1, a, b, b^3, a^2 b^2 \}, \,\,  \{ a^2 b^6, b^6, a^2 b^3, a^2 b^4, a \}, \,\, \{ 1, a^2, a b, b^2, b^6 \},\\[1mm]
\phantom{\{\,\, } \{ a^2 b, a b, b^5, a^2 b^2, a^2 b^4 \}, \,\, \{ 1, a^2 b, a^2 b^5, a b^6, a^2 b^6 \}, \,\, \{ a b^6, b, b^2, a^2 b^4, b^4 \}, \\[1mm]
\phantom{\{\,\, } \{ 1, a b^3, b^4, a^2 b^3, b^5 \}, \,\,  \{ a^2 b^5, b^3, a^2, a b^3, a^2 b^4 \}, \,\, \{ b, a^2, a^2 b, a^2 b^3, a b^5 \}, \\[1mm]
\phantom{\{\,\, } \{ a b^5, b^3, b^5, b^2, a^2 b^6 \}, \,\, \{ b, a b^2, b^5, b^6, a^2 b^5 \}, \,\,  \{ a^2, b^4, a^2 b^2, a^2 b^6, a b^2 \}, \\[1mm]
\phantom{\{\,\, } \{ b^2, a^2 b^2, a^2 b^3, a b^4, a^2 b^5 \}, \,\, \{ a b^4, b^4, b^6, a^2 b, b^3 \}, \,\, \{ 1, a b^2, a b^4, a^2 b^4, a b^5 \}, \\[1mm]
\phantom{\{\,\, } \{ a, a^2, a b^4, b^5, a b^6 \}, \,\, \{ a, a b, b^4, a b^5, a^2 b^5 \}, \,\, \{ a, b^2, a^2 b, a b^2, a b^3 \}, \\[1mm]
\phantom{\{\,\, } \{ b, a b, a b^3, a b^4, a^2 b^6 \}, \,\, \{ a b, a b^2, b^3, a^2 b^3, a b^6 \}, \,\, \{ a^2 b^2, a b^3, a b^5, b^6, a b^6 \} \,\, \}.
\end{array}$$
By Theorem~\ref{tmgrcube}, the design yields a $3$-cube $C_3$. We computed $|\Atop(C_3)|$ $=441$, and therefore
$C_3$ is not equivalent neither with $C_1$ nor with $C_2$. As before, by considering sub-$3$-cubes, we
see that the conclusion holds for all dimensions: $\C^n(21,5,1)$ contains a non-difference
group cube for every $n\ge 3$.
\end{example}

The inequivalence of the cubes $C_1$, $C_2$, and $C_3$
was shown by calculating their full autotopy groups. There is only one $(21,5,1)$ design~$\D_0$
up to isomorphism and duality, so these three cubes have the same slice invariant
$\{\{\D_0^{21}\}^3\}$ defined in Section~\ref{sec2}. For parameters $(16,6,2)$, there
are three inequivalent designs $\D_1$, $\D_2$, and $\D_3$~\cite{CD07} with full
automorphism groups of orders $|\Aut(\D_1)|=11520$, $|\Aut(\D_2)|=768$, and $|\Aut(\D_3)|=384$.
By Proposition~\ref{cordc}, the $27$ inequivalent difference $n$-cubes $(16,6,2)$ have
slice invariants $\{\{\D_i^{16}\}^e\}$ for $e={n\choose 2}(n-2)16^{n-3}$ and
$i\in\{1,2,3\}$. Our next goal is  to construct group cubes with different slice invariants,
which will then clearly be non-difference cubes.

The design $\D_1$ belongs to a family of symmetric designs with the symmetric
difference property (\emph{SDP designs}) studied by Kantor~\cite{WK75}. They
can be constructed from difference sets in elementary abelian groups
$G_m=\Z_2^{2m}$ of orders $4^m$ as follows. Any singleton, e.g.\ $D_1=\{00\}$,
is a $(4,1,0)$ difference set in the Klein four-group $G_1=\{00,01,10,11\}$.
Here we use additive notation and $00$ is the neutral element.
By the product construction of Mann~\cite[Lemma~6.3.1]{HBM65},
\begin{equation}\label{mannconstruction}
D_m=(D_{m-1}^c\times D_1)\cup (D_{m-1}\times D_1^c)
\end{equation}
is a $(4^m,2^{m-1}(2^m-1),2^{m-1}(2^{m-1}-1))$ difference set in
the group $G_m=G_{m-1}\times G_1$. The development of $D_m$ is an
SDP design denoted $\mathcal{S}^{-1}(2m)$ in~\cite{WK75}, and the design~$\D_1$
coincides with $\mathcal{S}^{-1}(4)$.

The other two $(16,6,2)$ designs $\D_2$ and $\D_3$ can be obtained
by switching submatrices in the incidence matrix of $\D_1$:
$$\left(\begin{array}{cccc}
1 & 1 & 0 & 0\\
1 & 1 & 0 & 0\\
0 & 0 & 1 & 1\\
0 & 0 & 1 & 1\\
\end{array}\right) \,\,\,\longrightarrow\,\,\,
\left(\begin{array}{cccc}
0 & 0 & 1 & 1\\
0 & 0 & 1 & 1\\
1 & 1 & 0 & 0\\
1 & 1 & 0 & 0\\
\end{array}\right).$$
Clearly any such switch leaves the row and column sums invariant. It can be
checked that dot products of rows and columns are also not changed. To be specific,
number the group elements of $G_2=\Z_2^4$ lexicographically:
$g_0=0000$, $g_1=0001$, $g_2=0010$, $g_3=0011,\ldots$, $g_{15}=1111$.
Let $D_2=\{g_1,g_2,g_3,g_4,g_8,g_{12}\}$ be the difference set for
$\D_1$. Denote the blocks $B_i=g_i+D_2$, so that $\D_1=\{B_0,\ldots,B_{15}\}$.
Then the blocks $B_0=D_2$ and $B_1=\{g_0,g_2,g_3,g_5,g_9,g_{13}\}$
contain the elements $g_2$, $g_3$ and do not contain the elements $g_{14}$, $g_{15}$.
On the other hand, the blocks $B_{12}=\{g_0, g_4, g_8, g_{13}, g_{14}, g_{15}\}$
and $B_{13}=\{g_1, g_5, g_9, g_{12}, g_{14}, g_{15}\}$ do not contain
$g_2$, $g_3$ and do contain $g_{14}$, $g_{15}$. By replacing these four
blocks by their symmetric difference with $S_1=\{g_2,g_3,g_{14},g_{15}\}$
and leaving the other blocks unchanged, the design $\D_1$ is transformed
into~$\D_2$. Denote the new blocks $B'_i=B_i\Delta S_1$ for $i=0,1,12,13$
and $B'_i=B_i$ for $i=2,\ldots,11,14,15$. The design $\D_3$ is obtained by
applying another switch to $\D_2$. Take $S_2=\{g_6,g_7,g_{14},g_{15}\}$ and
replace the blocks $B''_i=B'_i\Delta S_2$, $i=0,1,4,5$ leaving the other
blocks unchanged. 
The designs $\D_2=\{B'_0,\ldots,B'_{15}\}$ and $\D_3=\{B''_0,\ldots,B''_{15}\}$
are not SDP designs.

Now notice that the changed blocks $B'_i$ and $B''_i$ are also $(16,6,2)$
difference sets in the group $G_2$, equivalent with the original difference
set $D_2$. Therefore, $\D_2$ and $\D_3$ are $(16,6,2)$ designs with all blocks
being difference sets in $G_2$. The development of any one of these blocks is
isomorphic to $\D_1$, thus $\D_2$ and $\D_3$ are not developments
of their blocks. By applying Theorem~\ref{tmgrcube} to these designs, we get
group cubes of dimension $n\ge 3$ with slice invariants
$\{\{\D_1^{16}\}^{e_1},\{\D_2^{16}\}^{e_2}\}\}$ and $\{\{\D_1^{16}\}^{e_1},\{\D_3^{16}\}^{e_2}\}\}$,
respectively, for $e_1={n-1\choose 2}(n-2)16^{n-3}$ and $e_2=(n-1)(n-2)16^{n-3}$.

We now generalise this construction to obtain non-difference cubes
in $\C^n(4^m,2^{m-1}(2^m-1),2^{m-1}(2^{m-1}-1))$. Let $A$ be an incidence matrix
of a symmetric design with these parameters and $J$ the all-ones matrix of order~$4^m$.
Then $J-A$ is the incidence matrix of the complementary design with
parameters $(4^m,2^{m-1}(2^m+1),2^{m-1}(2^{m-1}+1))$. The block matrix
\begin{equation}\label{recconstruction}
\left(\begin{array}{cccc}
J-A & A & A & A\\
A & J-A & A & A\\
A & A & J-A & A\\
A & A & A & J-A\\
\end{array}\right)
\end{equation}
is an incidence matrix of a $(4^{m+1},2^{m}(2^{m+1}-1),2^{m}(2^{n}-1))$ design.
This recursive construction was described in~\cite{RB65} and~\cite{WK75} by using
$\{-1,1\}$-matrices and the Kronecker product.

When applied to the incidence matrix of $\D_1$, the construction gives the
series of SDP designs $\mathcal{S}^{-1}(2m)$.
Applied to $\D_2$ and $\D_3$, the construction gives two other series of
nonisomorphic designs $\mathcal{S}^{-1}_2(2m)$ and $\mathcal{S}^{-1}_3(2m)$ with
the same parameters. The blocks of these designs are $(4^m,2^{m-1}(2^m-1),2^{m-1}(2^{m-1}-1))$
difference sets in $G_m$, because the first block row of~\eqref{recconstruction} is
obtained by the product construction~\eqref{mannconstruction} of rows of~$A$ with the
singleton $\{00\}$, the second block row with the singleton $\{01\}$, etc. The developments
of these difference sets are all isomorphic to $\mathcal{S}^{-1}(2m)$, so $\mathcal{S}^{-1}_2(2m)$
and $\mathcal{S}^{-1}_3(2m)$ are not developments of their blocks. Theorem~\ref{tmgrcube}
gives rise to two group $n$-cubes in~$G_m$. Slices with~$i_1$ fixed are incidence matrices
of $\mathcal{S}^{-1}(2m)$, and slices with $i_1$ varying are incidence matrices of
$\mathcal{S}^{-1}_2(2m)$ or $\mathcal{S}^{-1}_3(2m)$. Therefore, these two $n$-cubes
are not difference cubes and the following theorem holds.

\begin{theorem}\label{nondifcubes}
For every $m\ge 2$ and $n\ge 3$, the set $\C^n(4^m,2^{m-1}(2^m-1),2^{m-1}(2^{m-1}-1))$
contains at least two inequivalent group cubes that are not difference cubes.
\end{theorem}

The total number of inequivalent group cubes in $\C^n(4^m,2^{m-1}(2^m-1),2^{m-1}(2^{m-1}-1))$
is much larger; see Proposition~\ref{gc16prop} for $n=3$ and $m=2$. Notice that the parameters
are of Menon type and therefore by exchanging $0\to -1$ the cubes are transformed into $n$-dimensional
proper Hadamard matrices. These matrices are \emph{totally regular}, meaning that any
slice is a regular Hadamard matrix. The designs $\D_1$, $\D_2$, $\D_3$ are transformed into
inequivalent Hadamard matrices $H_1$, $H_2$, $H_3$. Up to equivalence, there are only two
other Hadamard matrices $H_4$ and $H_5$ of order $16$~\cite{NSweb}. The series of designs
$\mathcal{S}^{-1}(2m)$, $\mathcal{S}^{-1}_2(2m)$, $\mathcal{S}^{-1}_3(2m)$ are transformed
into Hadamard matrices which are Kronecker products of $H_1$, $H_2$, $H_3$ with the cyclic
Hadamard matrix of order~$4$, and are also not equivalent. In conclusion, the
two non-difference cubes of Theorem~\ref{nondifcubes} can be transformed into
proper $n$-dimensional Hadamard matrices with inequivalent slices, which cannot be
obtained by \cite[Theorem~2.9]{dL90} and other known constructions.

\section{Small examples}\label{sec5}

\begin{table}[!b]
\begin{center}
\begin{tabular}{ccccccc}
\hline
ID & Structure & Nds & Ndc & $\dev D$ & Tds & Ngc\\
\hline
1 & $\Z_{16}$ & 0 & 0 & -- & 0 & 0 \\
2 & $\Z_4^2$ & 3 & 3 & $\D_1$ & 192 & 55 \\
3 & $(\Z_4 \times \Z_2) \rtimes \Z_2$ & 4 & 4 & $\D_1$ & 192 & 83 \\
4 & $\Z_4 \rtimes \Z_4$ & 3 & 3 & $\D_1$ & 192 & 81 \\
5 & $\Z_8 \times \Z_2$ & 2 & 2 & $\D_1$, $\D_2$ & 192 & 106 \\
6 & $\Z_8 \rtimes \Z_2$ & 2 & 2 & $\D_1$ & 64 & 34 \\
7 & $D_{16}$ & 0 & 0 & -- & 0 & 0 \\
8 & $QD_{16}$ & 2 & 2 & $\D_1$ & 128 & 50 \\
9 & $Q_{16}$ & 2 & 2 & $\D_1$ & 256 & 71 \\
10 & $\Z_4 \times \Z_2^2$ & 2 & 2 & $\D_1$ & 448 & 131 \\
11 & $\Z_2 \times D_8$ & 2 & 2 & $\D_1$ & 192 & 52 \\
12 & $\Z_2 \times Q_8$ & 2 & 2 & $\D_1$, $\D_3$ & 704 & 197 \\
13 & $(\Z_4 \times \Z_2) \rtimes \Z_2$ & 2 & 2 & $\D_1$, $\D_3$ & 320 & 77 \\
14 & $\Z_2^4$ & 1 & 1 & $\D_1$ & 448 & 9\\
\hline
\end{tabular}
\vskip 2mm
\caption{The group $3$-cubes of order $v=16$.}\label{gc16}
\end{center}
\end{table}

The smallest parameters of symmetric designs are $(7,3,1)$, $(11,5,2)$,
and $(13,4,1)$. There is a single design up to isomorphism and duality for each of
these parameters, coming from a difference set in the cyclic group $\Z_v$.
The next parameters are $(15,7,3)$ with four designs up to isomorphism
and duality, but only one coming from a difference set, namely $PG_2(3,2)$.
In all four cases the only designs with difference sets as blocks are
developments. Hence, there is a unique group cube in $\C^n(7,3,1)$, $\C^n(11,5,2)$,
$\C^n(13,4,1)$, and $\C^n(15,7,3)$, equivalent with the difference cube.

The fifth smallest parameters are $(16,6,2)$, and in this case there are many
non-difference group cubes. The results of an exhaustive computer search for
designs with $(16,6,2)$ difference sets as blocks are given in Table~\ref{gc16}. The first two columns
contain IDs of the groups of order~$16$ in the GAP library of small groups~\cite{GAP4} and a
description of their structure. The third and fourth column (Nds and Ndc) contain numbers of
inequivalent difference sets according to~\cite{DP19} and numbers of inequivalent difference cubes,
respectively. These numbers coincide by the computation from Section~\ref{sec3}. The fifth column contains designs
arising as developments of difference sets in each of the groups. The total number of difference sets,
including equivalent ones, is given in the sixth column (Tds). The last column (Ngc) contains numbers
of inequivalent non-difference $3$-cubes constructed by Theorem~\ref{tmgrcube}. For example, there are $58$
inequivalent group $3$-cubes in $\Z_4^2=\Z_4\times \Z_4$, three of which are difference cubes,
and $55$ of which are non-difference cubes. Here is a summary of the classification.

\begin{proposition}\label{gc16prop}
Up to equivalence, the set $\C^3(16,6,2)$ contains exactly $27$ difference cubes and
$946$ group cubes that are not difference cubes.
\end{proposition}

The most interesting examples come from groups with IDs $5$, $12$, and $13$, allowing
difference sets for inequivalent designs. For example, a design over $\Z_8\times \Z_2$
(group ID $5$) isomorphic to $\D_3$ can be constructed so that $8$ of its blocks
are difference sets with development~$\D_1$, and the other $8$ blocks are difference
sets with development~$\D_2$. It gives rise to a $3$-cube with slice invariant
$\{ \{ \D_3^{16} \}^2, \{ \D_1^{8}, \D_2^{8} \}^1 \}$ by Theorem~\ref{tmgrcube}. All
slice invariants that occurred are given in Table~\ref{c16inv}, along with the
corresponding numbers of difference cubes (Ndc) and non-difference group cubes (Ngc).

\begin{table}[t]
\begin{center}
\begin{tabular}{lccc}
\hline
Slice invariant & Ndc & Ngc & Nng \\
\hline
$\{ \{ \D_1^{16} \}^3 \}$ & 24 & 180 & 383 \\
$\{ \{ \D_2^{16} \}^3 \}$ & 1 & 20 & 32 \\
$\{ \{ \D_3^{16} \}^3 \}$ & 2 & 41 & 1 \\
$\{ \{ \D_1^{16} \}^2, \{ \D_2^{16} \}^1 \}$ & 0 & 15 & 392 \\
$\{ \{ \D_1^{16} \}^2, \{ \D_3^{16} \}^1 \}$ & 0 & 46 & 2 \\
$\{ \{ \D_2^{16} \}^2, \{ \D_1^{16} \}^1 \}$ & 0 & 284 & 444 \\
$\{ \{ \D_2^{16} \}^2, \{ \D_3^{16} \}^1 \}$ & 0 & 53 & 0 \\
$\{ \{ \D_3^{16} \}^2, \{ \D_1^{16} \}^1 \}$ & 0 & 189 & 77 \\
$\{ \{ \D_3^{16} \}^2, \{ \D_2^{16} \}^1 \}$ & 0 & 14 & 0 \\
$\{ \{ \D_1^{16} \}^2, \{ \D_1^{8}, \D_2^{8} \}^1 \}$ & 0 & 6 & 72 \\
$\{ \{ \D_1^{16} \}^2, \{ \D_1^{8}, \D_3^{8} \}^1 \}$ & 0 & 15 & 0 \\
$\{ \{ \D_2^{16} \}^2, \{ \D_1^{12}, \D_2^{4} \}^1 \}$ & 0 & 4 & 0 \\
$\{ \{ \D_2^{16} \}^2, \{ \D_1^{8}, \D_2^{8} \}^1 \}$ & 0 & 6 & 0 \\
$\{ \{ \D_2^{16} \}^2, \{ \D_1^{4}, \D_2^{12} \}^1 \}$ & 0 & 4 & 0 \\
$\{ \{ \D_2^{16} \}^2, \{ \D_1^{12}, \D_3^{4} \}^1 \}$ & 0 & 10 & 0 \\
$\{ \{ \D_2^{16} \}^2, \{ \D_1^{8}, \D_3^{8} \}^1 \}$ & 0 & 15 & 0 \\
$\{ \{ \D_2^{16} \}^2, \{ \D_1^{4}, \D_3^{12} \}^1 \}$ & 0 & 10 & 0 \\
$\{ \{ \D_3^{16} \}^2, \{ \D_1^{12}, \D_2^{4} \}^1 \}$ & 0 & 2 & 0 \\
$\{ \{ \D_3^{16} \}^2, \{ \D_1^{8}, \D_2^{8} \}^1 \}$ & 0 & 4 & 0 \\
$\{ \{ \D_3^{16} \}^2, \{ \D_1^{4}, \D_2^{12} \}^1 \}$ & 0 & 2 & 0 \\
$\{ \{ \D_3^{16} \}^2, \{ \D_1^{12}, \D_3^{4} \}^1 \}$ & 0 & 6 & 0 \\
$\{ \{ \D_3^{16} \}^2, \{ \D_1^{8}, \D_3^{8} \}^1 \}$ & 0 & 14 & 0 \\
$\{ \{ \D_3^{16} \}^2, \{ \D_1^{4}, \D_3^{12} \}^1 \}$ & 0 & 6 & 0 \\
$\{ \{ \D_1^{4}, \D_2^{12} \}^3 \}$ & 0 & 0 & 8 \\
$\{ \{ \D_1^{12}, \D_2^{4} \}^2, \{ \D_1^{4}, \D_2^{12} \}^1 \}$ & 0 & 0 & 12 \\
\hline
\end{tabular}
\vskip 2mm
\caption{Slice invariants of cubes in $\C^3(16,6,2)$.}\label{c16inv}
\end{center}
\end{table}

We have also found cubes of symmetric designs not equivalent with
any group cube, which we shall call \emph{non-group cubes}. Using our
package PAG~\cite{PAG}, we computed the full autotopy groups of available cubes
$C\in \C^3(16,6,2)$. We then chose subgroups $G\le \Atop(C)$ and used the Kramer-Mesner
method to construct all $3$-cubes of $(16,6,2)$ designs with $G$ as prescribed autotopy
group. Clearly we must always get the cube $C$ we started from, but often we
also get other inequivalent cubes, some of which are non-group cubes. Here is an
example.

\begin{example}\label{ngc}
Let $G=\langle \alpha_1, \alpha_2, \alpha_3 \rangle$ be the group
generated by the permutations
$$\begin{array}{rl}
\alpha_1= & ( 1,16)( 4, 5)( 6,11)( 7, 9)( 8,10)(14,15)\\
 & (17,28)(20,21)(22,27)(23,26)(24,25)(31,32)\\
 & (33,44)(34,37)(35,36)(38,39)(40,41)(47,48),\\[2mm]
\alpha_2= & ( 1,14, 2)( 3,16,15)( 4,13, 6)( 5,12,11)( 8, 9,10)\\
 & (17,20,29)(18,27,32)(19,22,31)(21,30,28)(23,24,25)\\
 & (33,47,46)(34,36,37)(38,40,42)(39,41,43)(44,48,45),\\[0mm]
\end{array}$$
$$\begin{array}{rl}
\alpha_3= & ( 1,13)( 2,11) ( 3, 6)( 7, 8)(12,16)(14,15)\\
 & (17,30,27,18)(19,28,29,22)(20,32,21,31)(23,25,24,26)\\
 & (33,43,38,46)(34,36,35,37)(39,45,44,42)(40,48,41,47).\\[0mm]
\end{array}$$
This is a group of order $384$ isomorphic to $\Z_2^6\rtimes S_3$ acting
on the sets $\P_1=\{1,\ldots,16\}$, $\P_2=\{17,\ldots,32\}$, and $\P_3=\{33,\ldots,48\}$.
The sets are point classes of a transversal design with blocks
being $G$-orbits generated by the subsets
\begin{align*}
\{ 1, 17, 33 \},\,\,\, \{ 1, 17, 40 \},\,\,\, \{ 1, 18, 33 \},\,\,\, \{ 1, 18, 34 \},\,\,\, \{ 1, 18, 42 \}, \\[1mm]
\{ 1, 23, 34 \},\,\,\, \{ 1, 23, 40 \},\,\,\, \{ 7, 17, 35 \},\,\,\, \{ 7, 17, 40 \},\,\,\, \{ 7, 23, 33 \}.
\end{align*}
This transversal design represents a $3$-cube of $(16,6,2)$ designs with slice invariant
$\{ \{ \D_1^{4}, \D_2^{12} \}^3 \}$. It is not equivalent with any group cube, because
parallel slices in all three directions comprise both designs $\D_1$ and $\D_2$.
A $3$-cube constructed from Theorem~\ref{tmgrcube} must have slices isomorphic to~$\D$
in two directions.
\end{example}

We found many non-group cubes of $(16,6,2)$ designs by this method.

\begin{proposition}\label{ngc16prop}
The set $\C^3(16,6,2)$ contains at least $1423$ inequivalent non-group cubes.
\end{proposition}

The last column of Table~\ref{c16inv} (Nng) contains their distribution by slice invariants.
The constructed examples are available on our web page:
\begin{center}
\url{https://web.math.pmf.unizg.hr/~krcko/results/cubes.html}
\end{center}
The total number of non-group cubes in $\C^3(16,6,2)$ is probably much larger. We
have attempted construction only for groups of order $|G|\ge 512$.

We did not find any other non-group cubes for smaller parameters $(v,k,\lambda)$,
and for larger parameters the Kramer-Mesner approach was too inefficient.
However, we did find larger non-difference group cubes by the construction of
Theorem~\ref{tmgrcube}. A complete classification was possible for $(21,5,1)$.
The non-difference group cube $C_3$ of Example~\ref{pg24} is unique.

\begin{proposition}\label{gc21prop}
The set $\C^3(21,5,1)$ contains exactly three inequivalent group
cubes, two of which are difference cubes.
\end{proposition}

\begin{table}[!b]
\begin{center}
\begin{tabular}{cccc}
\hline
Parameters & Nds & Ndc & Ngc \\
\hline
$(21,5,1)$ & $2$ & $2$ & $1$ \\
$(27,13,6)$ & $3$ & $2$ & $\ge 7$ \\
$(36,15,6)$ & $35$ & $35$ & $\ge 373$ \\
$(45,12,3)$ & $2$ & $2$ & $\ge 6$ \\
$(63,31,15)$ & $6$ & $6$ & $\ge 9$ \\
$(64,28,12)$ & $330159$ & $<330159$ & $\ge 1$ \\
$(96,20,4)$ & $2627$ & $1806$ & $\ge 1$ \\
\hline
\end{tabular}
\vskip 2mm
\caption{Difference and group cubes for $n=3$.}\label{tabndgc}
\end{center}
\end{table}

For parameters $(27,13,6)$ we did an incomplete search and found $7$ other group cubes
besides the two difference cubes mentioned at the end of Section~\ref{sec3}. For even
larger parameters we searched for designs with difference sets as blocks and a prescribed
automorphism group. The parameters for which we found non-difference group cubes are given
in Table~\ref{tabndgc}. An on-line version of the table with links to files containing
the cubes is available on the above-mentioned web page. The column Nds contains numbers
of inequivalent difference sets according to~\cite{DP19}. For parameters $(64,28,12)$
some of the corresponding difference $3$-cubes are equivalent, but we did not perform a
complete enumeration.

\end{document}